\documentclass[12pt,bezier]{amsart}
\usepackage{amsfonts,latexsym,amscd}
\newtheorem{theorem}{Theorem}[section]

\newtheorem{proposition}[theorem]{Proposition}
\newtheorem{corollary}[theorem]{Corollary}
\newtheorem{lemma}[theorem]{Lemma}

\newcommand{\thm}[1]{Theorem \ref{#1}}
\newcommand{\eqn}[1]{Equation \ref{#1}}
\newcommand{\prop}[1]{Proposition \ref{#1}}
\newcommand{\lma}[1]{Lemma \ref{#1}}

\newcommand{\cor}[1]{Corollary \ref{#1}}

\newcommand{\rarrow}{\rightarrow}
\newcommand{\fd}{\operatorname{fd}}

\newcommand{\spec}{\operatorname{Spec}}

\newcommand{\ho}{\operatorname{Ho}}

\newcommand{\chr}{\operatorname{char}}
\newcommand{\depth}{\operatorname{depth}}
\newcommand{\tor}{\operatorname{Tor}}
\newcommand{\wgt}{\operatorname{wt}}
\newcommand{\aqdim}{\operatorname{AQ-dim}}
\newcommand{\la}{\operatorname{\!\langle\!}}
\newcommand{\ra}{\operatorname{\!\rangle}}

\newcommand{\F}{{\mathbb F}}

\newcommand{\fq}{\mathfrak{q}}

\newcommand{\fm}{\mathfrak{m}}

\newcommand{\calA}{{\mathcal A}}
\newcommand{\Ap}{{\mathcal A}_{\ell}}

\newcommand{\calS}{{\mathcal S}}
\newcommand{\calN}{{\mathcal N}}

\newcommand{\calV}{{\mathcal V}}

\begin{document}

\title[On Simplicial Algebras with Finite Andr\'e-Quillen Homology]
{On Simplicial Commutative Algebras with Finite Andr\'e-Quillen 
Homology}
\author{James M. Turner}
\address{Department of Mathematics\\
Calvin College\\ 
3201 Burton Street, S.E.\\
Grand Rapids, MI 49546}
\email{jturner@calvin.edu}
\thanks{Partially supported by National Science Foundation (USA) grant 
DMS-0206647 and a Calvin Research Fellowship. He thanks the Lord for 
making his work possible.}
\date{\today}

\begin{abstract}
In \cite{Qui2, Avr} a conjecture was posed to the effect that if $R 
\to A$ is a homomorphism of Noetherian rings then the 
Andr\'e-Quillen homology on the category of $A$-modules satisfies:
$D_{s}(A|R;-) = 0$ for $s\gg 0$ implies $D_{s}(A|R;-) = 0$ for $s\geq 3$. In 
\cite{Tur3}, an extended version of this conjecture was considered for 
which $A$ is a simplicial commutative $R$-algebra with Noetherian 
homotopy such that $\chr (\pi_{0}A) \neq 0$. In addition, a homotopy characterization 
of such algebras was described. The main goal of this paper is to develop 
a strategy for establishing this extended conjecture and
provide a complete proof when $R$ is Cohen-Macaulay
of characteristic 2.
\end{abstract}

\maketitle

\section*{Overview}
In \cite{Qui2}, D. Quillen presented his viewpoint on the homology of 
algebras which extended, in the commutative case, the work of 
Lichtenbaum and Schlessinger and gave M. Andr\'e's notion of 
homology. Furthermore, he observed that strong vanishing of this
Andr\'e-Quillen homology for finite type algebras held only when such algebras 
possessed the complete intersection property and conjectured that a 
weaker type of vanishing also characterized such algebras. In 
\cite{Avr}, L. Avramov clarified and extended Quillen's conjectures in 
the following manner. Let $f: R \to A$ be a homomorphism of Noetherian 
rings [Note: unless otherwise noted, all rings and algebras from this 
point on are commutative with unit]. Then $f$ is a {\it locally complete 
intersection} provided for each $\fq\in \spec S$ the semi-completion 
$R_{\wp\cap R}\to \hat{A_{\wp}}$ suitably factors through a 
surjection with kernel being generated by a regular sequence (see below 
for more details).

\bigskip

\noindent {\bf Quillen's Conjecture:} (see \cite{Avr, Qui2}) {\it Let $R \to A$ 
be a homomorphism of Noetherian rings such that the Andr\'e-Quillen 
homology satisfies $D_{s}(A|R;-) = 0$ (as functors of $A$-modules) for 
$s\gg 0$. Then
\begin{enumerate}
	\item $D_{s}(A|R;-) = 0$ for $s\geq 3$;
	\item if $\fd_{R}A < \infty$ then $R \to A$ is a locally complete 
	intersection (and, hence, $D_{s}(A|R;-) = 0$ for $s\geq 2$).
\end{enumerate}}

\bigskip

\noindent Part 2 of this conjecture was proved by Avramov in \cite{Avr}. Part 1 was 
proved by Avramov and S. Iyengar for algebra retracts in \cite{AI1}. 

Following ideas of Haynes Miller, an alternate approach to proving 
this conjecture was taken in \cite{Tur1, Tur2} when $R$ is a field 
by viewing it as a special case of an algebraic version of a theorem 
of J.P. Serre \cite{Serre}. Following this line of thinking, in 
\cite{Tur2,Tur3} the more general consideration of Noetherian algebras was extended to 
{\it simplicial commutative algebras with Noetherian homotopy}, that 
is, simplicial commutative algebras $A$ such that $\pi_{0}A$ is 
Noetherian and $\pi_{*}A$ is a finite graded $\pi_{0}A$-module. In 
using Andr\'e-Quillen homology to analyse such, we can use the type of 
tools first clarified by Andr\'e and Quillen: flat base change, 
transitivity sequence, localization etc. Cf. \cite{And, Qui2, Tur3}. 
A particularly useful method for analysing simplicial commutative 
algebras in our present context through homology is the following 
generalization of the main result in \cite{AFH}, proved in \cite{Tur3}.
For each $\wp \in \spec (\pi_{0}A)$ the simplicial commutative 
algebra $A^{\prime} = A\otimes^{{\bf L}}_{\pi_{0}A}\widehat{(\pi_{0}A)_{\wp}}$ 
there is a complete local ring $R^{\prime}$ and a homotopy commutative diagram

$$
\begin{array}{ccc}
R & \stackrel{\eta}{\longrightarrow} &
A \\[1mm]
\phi \downarrow \hspace*{10pt}
&&
\hspace*{10pt} \downarrow \psi \\[1mm]
R^{\prime} & \stackrel{\eta^{\prime}}{\longrightarrow} & A^{\prime}
\end{array}
$$
with the following properties:
	 \begin{enumerate}
		 
		 \item $\phi$ is a flat map and its closed fibre $R^{\prime}/\wp 
		 R^{\prime}$ is weakly regular;
		 
		 \item $\psi$ is a $D_{*}(-|R;k(\wp))$-isomomorphism;
		 		 
		 \item $\eta^{\prime}$ induces a surjection $\eta^{\prime}_{*}: R^{\prime} \to 
		 (\pi_{0}A^{\prime}, k(\wp))$ of local rings;
		 
		 \item $\fd_{R} (\pi_{*} A)$ finite implies that $\fd_{R^{\prime}} 
		 (\pi_{*} A^{\prime})$ is finite
		 
	 \end{enumerate}	 
We call such a diagram a {\it homotopy factorization} of $A$. We can use 
such factorizations to extend the notion of locally complete intersection to 
simplicial commutative $R$-algebras with Noetherian homotopy. 
Specifically, we call such $A$ a {\it a locally homotopy 
n-intersection}, $n$ a natural number, provided for each $\wp\in 
\spec (\pi_{0}A)$ there is a factorization such that the {\it 
connected component at $\wp$} satisfies
$$
A(\wp) := A^{\prime}\otimes^{\bf L}_{R^{\prime}}k(\wp) \simeq S_{k(\wp)}(W)
$$
with $W$ a connected simplicial $k(\wp)$-module satisfying 
$\pi_{s}W = 0$ for $s>n$. Here and throughout $S_{k(\wp)}(-)$ denotes the 
free commutative $k(\wp)$-algebra functor.

We can now state, inspired by Serre's theorem \cite{Serre}, our simplicial version of 
Quillen's conjecture:

\bigskip

\noindent {\bf Vanishing Conjecture:} {\it Let $R$ be a Noetherian 
ring and let $A$ be a simplicial commutative $R$-algebra with finite 
Noetherian homotopy and $\chr (\pi_{0}A) \neq 0$ such that the 
Andr\'e-Quillen homology satisfies $D_{s}(A|R;-) = 0$ 
(as functors of $\pi_{0}A$-modules) for $s\gg 0$. Then
\begin{enumerate}
	\item $A$ is a locally homotopy 2-intersection;
	\item if $\fd_{R}\pi_{*}A < \infty$ then $A$ is a locally homotopy 
	1-intersection.
\end{enumerate}	
}

\bigskip
 
Part 2 of the Vanishing Conjecture was proved in \cite{Tur2,Tur3}. The 
goal of the first part of this paper is to outline a strategy for 
giving a proof for the whole Vanishing Conjecture. The strategy 
involves formulating a more local version of the Algebraic Serre 
Theorem proved in \cite{Tur2} used to prove part 2 of the Vanishing 
Conjecture. Specifically, we will analyse the behavior of homotopy 
operations on each $A(\wp)$, particular the divided $p^{th}$-powers 
and the {\it Andr\'e operation} (so named because of the role it 
played in \cite{And3} which motivated the direction of this paper). 
In the first section, we will formulate a Nilpotence and 
Non-nilpotence Conjecture regarding the action of these operations 
which, when coupled together, imply the Vanishing Conjecture. The 
second section will then focus on proving the Nilpotence Conjecture 
at the prime 2. Finally, in the third section we will establish what 
will hopefully be our first step toward proving the Vanishing 
Conjecture when $\chr (\pi_{0}A) = 2$. Specifically, we will 
establish our:

\bigskip

\noindent {\bf Main Theorem:} {\it Let $A$ be a simplicial commutative 
$R$-algebra with finite Noetherian homotopy such that $R$ is Cohen-Macaulay 
of characteristic 2. Then $D_{s}(A|R;-) = 0$ (as a functor of 
$\pi_{0}A$-modules) for $s\gg 0$ if and only if $A$ is a locally 
homotopy 2-intersection.}

\bigskip 

As an immediate consequence, we obtain:

\bigskip

\noindent {\bf Corollary.} {\it
	Let $R \to A$ be a homomorphism of Noetherian rings of characteristic 
	2 such that $R$ is Cohen-Macaulay. Then $D_{s}(A|R;-) = 0$ for $s 
	\gg 0$ implies $D_{s}(A|R;-) = 0$ for $s \geq 3$.}

\bigskip

\noindent {\bf Acknowledgements.} The author would like to thank
Lucho Avramov for educating him on Cohen-Macaulay
rings and for comments and criticisms on an earlier draft of this
paper. He would also like to thank Paul Goerss for several discussions 
on homotopy operations as well as for many other helpful comments.

\section{Nilpotence Conjectures}

In this section we reformulate the Vanishing Conjecture in terms of a 
two part Nilpotence Conjecture which shifts the burden for global vanishing 
of Andr\'e-Quillen homology to local vanishing of operations acting on 
the homotopy of components. We will first need a weaker notion of 
homotopy factorization in order to tighten our grip on how on how 
information from the homotopy of our simplicial algebra is transferred 
to the homotopy of its components. We will, throughout this section, be 
assuming basic properties of Andr\'e-Quillen homology, refering the 
reader to \cite{Tur3} for details. 

\subsection{Weak homotopy factorizations}

In the next subsection, we will recall that the conclusions of the 
Vanishing Conjectures are equivalent to certain strong global 
vanishing properties of Andr\'e-Quillen homology. Our goal at present 
is to modify the notion of homotopy factorizations which suitably 
preserves the Andr\'e-Quillen homology but puts a tighter control on 
the local ring $R^{\prime}$.

Let $A$ be a simplicial commutative $R$-algebra and denote $\pi_{0}A$ by 
$\Lambda$. We may assume that $A$ is a simplicial commutative 
$\Lambda$-algebra. Cf. \cite[Theorem A]{Tur3}.
Fix $\wp \in \spec \Lambda$ and let $\widehat{(-)}$ denote the 
completion functor on $R$-modules at $\wp$. Define the homotopy 
connected simplicial supplemented $\widehat \Lambda$-algebra $A^{\prime}$ by
$$
A^{\prime} = A\otimes^{{\bf L}}_{\Lambda}\widehat \Lambda.
$$

\begin{proposition}\label{final}
	Suppose $A$ is a simplicial commutative $R$-algebra with $R$ a Noetherian ring. 
	Then there exists a (complete) local ring $(R^{\prime\prime},\fm)$, a 
	simplicial commutative $R^{\prime\prime}$-algebra $A^{\prime\prime}$, 
	and a homotopy commutative diagram 
\begin{equation}\label{wfac}
\begin{array}{ccc}
R & \stackrel{\eta}{\longrightarrow} &
A \\[1mm]
\phi \downarrow \hspace*{10pt}
&&
\hspace*{10pt} \downarrow \psi \\[1mm]
R^{\prime\prime} & \stackrel{\eta^{\prime\prime}}{\longrightarrow} & 
A^{\prime\prime}
\end{array}
\end{equation}
     with the following properties:
	 \begin{enumerate}
		 
		 \item $\phi$ is a complete intersection at $\fm$;
		 
		 \item $\depth (\fm) = 0$;
		 
		 \item $D_{\geq 2}(A|R;k(\wp))\cong D_{\geq 2}(A^{\prime\prime}|R^{\prime\prime};k(\wp))$;
		 
		 \item $\eta^{\prime\prime}$ induces a surjection of local rings
		 $\eta^{\prime\prime}_{*}: R^{\prime\prime} \to \pi_{0}A^{\prime\prime}$;
		 
		 \item If $A$ has finite Noetherian homotopy then $A^{\prime\prime}$ 
		 has finite Noetherian homotopy.
		 		 
	\end{enumerate}	 

\end{proposition}	

\noindent {\it Proof:} Choose a homotopy factorization of $A$ over 
$\wp$
$$
\begin{array}{ccc}
R & \stackrel{\eta}{\longrightarrow} &
A \\[1mm]
\phi \downarrow \hspace*{10pt}
&&
\hspace*{10pt} \downarrow \psi \\[1mm]
R^{\prime} & \stackrel{\eta^{\prime}}{\longrightarrow} & A^{\prime}
\end{array}
$$
which exists by \cite[(2.8)]{Tur3}. 

Next, let $\fq$ be the maximal ideal 
of $R^{\prime}$. Let 
$x_{1},\ldots,x_{r}$ be a maximal $R^{\prime}$-subsequence of a minimal
generating set for $\fq$. We define 
$$R^{\prime\prime} = R^{\prime}/(x_{1},\ldots,x_{r}).$$
Then $\fm = \fq/(x_{1},\ldots,x_{r})\fq$ has depth 0 since it 
contains only zero divisors. Furthermore, the composite $R_{\wp}\to R^{\prime} 
\to R^{\prime\prime}$ is a complete intersection at $\fm$ by definition. 
Cf. \cite{Avr}. 

Now, let $A^{\prime\prime} = A^{\prime}\otimes_{R^{\prime}}^{{\bf 
L}}R^{\prime\prime}$. Then
$$
D_{\geq 2}(A|R;k(\wp)) \cong D_{\geq 2}(A^{\prime}|R;k(\wp))
\cong D_{\geq 2}(A^{\prime}|R^{\prime};k(\wp)) \cong 
D_{\geq 2}(A^{\prime}\otimes_{R^{\prime}}^{{\bf L}}R^{\prime\prime}|R^{\prime\prime};k(\wp))
$$
which follows from the properties of homotopy factorizations, the 
transitivity sequence, and flat base change \cite[(2.4)]{Tur3}. 
Applying $\pi_{0}$ to the map $R^{\prime\prime} 
\to A^{\prime}\otimes_{R^{\prime}}^{{\bf L}}R^{\prime\prime}$ gives 
the map $R^{\prime\prime} \cong 
R^{\prime}\otimes_{R^{\prime}}R^{\prime\prime} \to 
\pi_{0}(A^{\prime})\otimes_{R^{\prime}}R^{\prime\prime}$ which is 
a surjection. Thus $R^{\prime\prime} \to \pi_{0}A^{\prime\prime}$ is a 
surjection. 

Finally, if $A$ has finite Noetherian homotopy then so does 
$A^{\prime}$ (since $\pi_{*}A^{\prime} \cong \widehat{\pi_{*}A}$). By 
\cite[\S II.6]{Qui1}, there is a Kunneth spectral sequence
$$
E^{2}_{s,t} = 
\tor_{s}^{R^{\prime}}(\pi_{t}A^{\prime},R^{\prime\prime}) 
\Longrightarrow \pi_{s+t}A^{\prime\prime}.
$$
Since $R^{\prime} \to R^{\prime\prime}$ is a complete intersection, 
$\fd_{R^{\prime}}R^{\prime\prime} < \infty$. Thus 
$\pi_{*}A^{\prime\prime}$ will be a finite module over 
$\pi_{0}A^{\prime\prime} \cong 
\widehat\Lambda\otimes_{R^{\prime}}R^{\prime\prime}$.
\hfill $\Box$

\bigskip

We will call a diagram (\ref{wfac}) satisfying the conditions (1) - 
(5) above a {\it weak homotopy factorization} for $A$. 

\subsection{Brief review of the homotopy of simplicial commutative algebras 
over a field}

Let $A$ be a simplicial commutative $\ell$-algebra where $\ell$ is a 
field. In this section we review some basic facts about the homotopy 
groups of such objects, computed as the homotopy groups of simplicial 
$\ell$-modules.

Let $\calA_{\ell}$ be the category of {\it supplemented 
$\ell$-algebras}, i.e. commutative $\ell$-algebras augmented over 
$\ell$. Let $s\Ap$ be the category of simplicial objects over 
$\calA_{\ell}$. Then for $A\in s\Ap$ and $n\geq 0$ we have a natural 
isomorphism
$$
\pi_{n}A \cong [S_{\ell}(n), A]_{\ho (s\Ap)}
$$
where $S_{\ell}(n) = S_{\ell}(K(n))$, $K(n)$ the simplicial 
$\ell$-module satisfying $\pi_{*}K(n) \cong \ell$ concentrated in 
degree $n$. We will use this relation to determine the natural primary 
algebra structure on $\pi_{*}A$. 

Given integers $r_{1},\ldots,r_{m},t_{1},\ldots,t_{n} \neq 0$ an 
{\it multioperation} of degree \\ $(r_{1},\ldots,r_{m};t_{1},\ldots,t_{n})$ 
is a natural map
$$
\theta : \pi_{r_{1}}\times \ldots \times \pi_{r_{m}} \to 
\pi_{t_{1}}\times \ldots\times \pi_{t_{n}}
$$
of functors on $s\Ap$. Let 
$\mbox{Nat}_{r_{1},\ldots,r_{m};t_{1},\ldots,t_{n}}$ be the 
set of multioperations of degree \\
$(r_{1},\ldots,r_{m};t_{1},\ldots,t_{n})$. It is straightforward to 
show that 
$$
\mbox{Nat}_{r_{1},\ldots,r_{m};t_{1},\ldots,t_{n}} \cong 
\mbox{Nat}_{r_{1},\ldots,r_{m};t_{1}}\times\ldots\times 
\mbox{Nat}_{r_{1},\ldots,r_{m};t_{n}}.
$$

Now, we define 
\begin{equation}\label{natmap}
f: \mbox{Nat}_{r_{1},\ldots,r_{m};t} \to 
\pi_{t}(S_{\ell}(r_{1})\otimes_{\ell}\ldots\otimes_{\ell} S_{\ell}(r_{m}))
\end{equation}
as follows. Let $\calN = \mbox{Nat}_{r_{1},\ldots,r_{m};t}$
and let $X = S_{\ell}(r_{1})\otimes_{\ell}\ldots\otimes_{\ell} S_{\ell}(r_{m}).$
For each $1\leq j \leq m$, let $\iota_{j}\in \pi_{r_{j}}X$ be the 
homotopy class of the inclusion $S_{\ell}(r_{j}) \to X$. Given 
$\theta \in \calN$ there is an induced map
$$
\theta_{X} : \pi_{r_{1}}X \times \ldots \times \pi_{r_{m}}X \to \pi_{t}X.
$$
Thus we can define $f : \calN \to \pi_{t}X$
by 
\begin{equation}\label{natmapdef}
f(\theta) = \theta_{X}(\iota_{1},\ldots,\iota_{m}).
\end{equation}

\begin{proposition}\label{natprop}
	$\mbox{Nat}_{r_{1},\ldots,r_{m};t} \cong 
	\pi_{t}(S_{\ell}(r_{1})\otimes_{\ell}\ldots\otimes_{\ell} S_{\ell}(r_{m})).$
\end{proposition}

\noindent {\it Proof:} Since we have
$$
\pi_{r_{1}} \times \ldots \times \pi_{r_{m}} \cong 
[S_{\ell}(r_{1})\otimes_{\ell}\ldots\otimes_{\ell} 
S_{\ell}(r_{m}),-]_{\ho (s\Ap)}
$$
the result follows from Yoneda's lemma \cite{Mac}. \hfill $\Box$

\bigskip

\noindent {\bf Note:} 
\begin{enumerate}
\item There is an obvious map 
$$
\mbox{Nat}_{r_{1},\ldots,r_{m};t}\times\mbox{Nat}_{t;q} \to 
\mbox{Nat}_{r_{1},\ldots,r_{m};q}
$$
induced by composition. 
\item $\mbox{Nat}$ is naturally an 
$\ell$-module and $f$ is naturally a linear map.
\end{enumerate}

We now can address the issue of understanding possible relations among  
multioperations.  

\begin{corollary}\label{natrel}
	For $\theta \in \mbox{Nat}_{r_{1},\ldots,r_{m};t}$ then any 
	expression for $\theta$ in $\mbox{Nat}_{r_{1},\ldots,r_{m};q}$ is 
	determined by 
	$f(\theta)\in \pi_{t}(S_{\ell}(r_{1})\otimes_{\ell}\ldots\otimes_{\ell} 
	S_{\ell}(r_{m}))$.
	Furthermore, if $\psi\in \mbox{Nat}_{t;q}$ then $f(\psi\circ 
	\theta) = f(\psi)\circ f(\theta)$, as composites of their homotopy 
	representitives, in 
	$\pi_{t}(S_{\ell}(r_{1})\otimes_{\ell}\ldots\otimes_{\ell} S_{\ell}(r_{m}))$.
\end{corollary}	

\noindent {\it Proof:} This again follows from Yoneda's lemma \cite{Mac}.
\hfill $\Box$

\bigskip

Now, we are in a position to determine the full natural primary 
structure for homotopy in $s\Ap$. First, recall that for any field 
$\F$ we have
\begin{equation}\label{symid}
S_{\F}(V\oplus W) \cong S_{\F}(V)\otimes S_{\F}(W).
\end{equation}
Next, let $k = {\mathbb Q}$ if $\chr (\ell) = 0$ and let $k = \F_{p}$ 
if $\chr (\ell) = p \neq 0$. We seek a natural map of $\ell$-algebras
$$
\phi_{V} : S_{\ell}(V\otimes_{k}\ell) \to S_{k}(V)\otimes_{k}\ell
$$
where $V$ is a $k$-module. This can be defined as the adjunction of 
the inclusion $V\otimes_{k}\ell \to I(S_{k}(V)\otimes_{k}\ell)$ 
(here $I: \calA_{\ell} \to \calV_{\ell}$ is the augmentation ideal 
functor).

\begin{proposition}\label{symprop}
	The natural map $\phi : S_{\ell}((-)\otimes_{k}\ell) \to 
	S_{k}(-)\otimes_{k}\ell$ is an isomorphism of functors from 
	$k$-modules to $\Ap$.
\end{proposition}	

\noindent {\it Proof:} By the identity (\ref{symid}) and naturality, it is 
enough to provide a proof for one dimensional $V$, i.e. for $V \cong k\la 
x\ra$. Then $\phi_{V} : \ell [x] \to k [x]\otimes_{k}\ell$ is 
determined algebraically by the value $\phi_{V}(x) = x\otimes_{k}1$. This is clearly 
an isomorphism. 
\hfill $\Box$

\bigskip

\noindent {\bf Note:} This and other similar results can also be 
shown to follow from the faithful flatness of the functor 
$(-)\otimes_{k}\ell$.

\bigskip

\begin{corollary}\label{homsymprop}
	For $V \in s\calV_{k}$ there is a natural isomorphism 
	$$\pi_{*}(S_{\ell}(V\otimes_{k}\ell)) \cong 
	\pi_{*}(S_{k}(V))\otimes_{k}\ell.$$ As a consequence all natural 
	primary homotopy operations for simplicial supplemented 
	$\ell$-algebras and their relations are determined by 
	$\pi_{*}S_{k}(n)$ for all $n\in {\bf N}$.
\end{corollary}	

\noindent {\it Proof:} The first statement follows from 
\prop{symprop} and the faithful flatness of $(-)\otimes_{k}\ell$. The second 
statement follows additionally from \cor{natrel} and the Kunneth 
theorem. Recall that $S_{\ell}(n) \cong S_{\ell}(K_{\ell}(n))$ and we 
can take $K_{\ell}(n) = \ell\la S^{n} \ra \cong k\la S^{n} 
\ra\otimes_{k}\ell$, where $S^{n}$ is a choice of simplicial set model 
for the n-sphere. 
\hfill $\Box$

\bigskip

\noindent {\bf Note:} The computation of $\pi_{*}S_{{\mathbb Q}}(n)$ can 
be traced back to at least as early as \cite{Car}. The computation 
of $\pi_{*}S_{\F_{p}}(n)$ can be found in \cite{Bou, Car}, for
general $p$, and in \cite{Dwy} for $p = 2$. We will review the 
results of \cite{Dwy} in the next section.

\bigskip

For non-zero characteristics, we will be interested in two 
particular operations. Specifically, for $A \in s\Ap$, $\pi_{*}A$ is 
naturally a divided power algebra. Therefore, there is a {\it 
divided $p^{th}$-power operation}
$$
\gamma_{p} : \pi_{n}A \to \pi_{pn}A.
$$
Cartan, Bousfield, and Dwyer also construct an operation 
$$
\vartheta : \pi_{n}A \to \pi_{(p-1)n+1}A
$$
which we call the {\it Andr\'e operation} because of the role it 
played in \cite{And3} where M. Andr\'e's showed that Gulliksen's result about 
the equality of deviations with simplicial dimensions for rational local 
rings \cite{Gul} cannot be extended to the primary case.
In the notion of \cite[8.8]{Bou} and \cite{Dwy}, 
\begin{equation}
  \vartheta =
  \begin{cases}
    \delta_{n-1} & \quad p = 2,\\
	\nu_{(n-1)/2} & \quad p > 2.
\end{cases}
\label{andopdef}
\end{equation}
A useful basic relation between the two operations is
\begin{equation}\label{anddivrel}
	\vartheta\gamma_{p} = 0
\end{equation}

\subsection{Nilpotency conjectures and consequences}

We now are in a position to address the Vanishing Conjecture and
reformulate it in terms of local conditions on homotopy groups. To 
begin, we need the following

\begin{lemma}\label{lowdim}
	Let $W \in s\calV_{\ell}$, with $\chr (\ell) > 0$, and let 
	$n \in {\mathbb N}$ be so that $\pi_{j}W \neq 0$ implies $n\geq j \geq 1$. 
	Then
	\begin{enumerate}
		\item $\gamma_{p} = 0$ on $\pi_{*}S_{\ell}(W)$ provided $n = 1$;
		\item $\vartheta = 0$ on $\pi_{*}S_{\ell}(W)$ provided $n = 2$.
	\end{enumerate}
\end{lemma}

\noindent {\it Proof:} By \cor{homsymprop}, it is enough to provide a 
proof for $\ell = \F_{p}$. For $n = 1$, $\pi_{*}S_{\ell}(W)$ is a 
free exterior algebra generated by $\pi_{1}W$, which has trivial 
$\gamma_{p}$-action. For $n = 2$, $\pi_{*}S_{\ell}(W)$ is a free 
divided power algebra generated by $\pi_{*}W$. Cf. \cite{Car}. Thus
$\pi_{*}S_{\ell}(W)$ has trivial $\vartheta$-action by relation 
(\ref{anddivrel}). \hfill $\Box$ 

\bigskip

Given $A \in s\Ap$ with $\chr (\ell) > 0$, we call $A$ 
\begin{enumerate}
	\item {\it $\Gamma$-nilpotent} provided $\gamma_{p}^{s}(x) = 0$ for 
	$s\gg 0$ for all $x \in \pi_{*}A$ and
	\item {\it Andr\'e nilpotent} provided $\vartheta^{s}(x) = 0$ for 
	$s\gg 0$ for all $x \in \pi_{*}A$.
\end{enumerate}
Next, given a Noetherian ring $R$ and a simplicial commutative 
$R$-algebra $A$ with Noetherian homotopy. For $\wp \in \spec(\pi_{0}A)$ 
with $\chr (k(\wp)) > 0$, we call $A$
\begin{enumerate}
	\item {\it $\Gamma$-nilpotent at $\wp$} provided there is a homotopy 
	factorization at $\wp$ such that $A(\wp)$
	is $\Gamma$-nilpotent over $k(\wp)$, and
	\item {\it Andr\'e nilpotent at $\wp$} provided there is a weak homotopy 
	factorization at $\wp$ such that 
	$A^{\prime\prime}\otimes_{R^{\prime\prime}}k(\wp)$
	is Andr\'e nilpotent over $k(\wp)$.
\end{enumerate}

\begin{proposition}
	Let $A$ be a simplicial commutative $R$-algebra with Noetherian 
	homotopy and $\wp \in \spec (\pi_{0}A)$ such that $\chr (k(\wp)) > 0$. Then 
	\begin{enumerate}
		\item $A$ is $\Gamma$-nilpotent at $\wp$ provided $A$ is a  
		homotopy 1-intersection at $\wp$; 
		\item $A$ is Andr\'e-nilpotent at $\wp$ provided $A$ is a  
		homotopy 2-intersection at $\wp$.
	\end{enumerate}
\end{proposition}

\noindent {\it Proof:} Both follow from the definitions and 
\lma{lowdim}. \hfill $\Box$

\bigskip

We now can state our two nilpotence-type conjectures.
\bigskip

\noindent {\bf Nilpotence Conjecture:} {\it Let $A$ be a simplicial 
commutative $R$-algebra with finite Noetherian homotopy. Let $\wp \in \spec 
(\pi_{0}A)$ be such that $\chr (k(\wp)) > 0$ and $D_{s}(A|R;k(\wp)) = 0$ 
for $s\gg 0$. Then:
\begin{enumerate}
	\item $A$ is $\Gamma$-nilpotent at $\wp$ if and only if $A$ is a 
	homotopy 1-intersection at $\wp$;
	\item $A$ is Andr\'e nilpotent at $\wp$ if and only if $A$ is a  
	homotopy 2-intersection at $\wp$.
\end{enumerate}}

\bigskip

\noindent {\bf Non-Nilpotence Conjecture:} {\it Let $A$ be a simplicial 
commutative $R$-algebra with finite Noetherian homotopy. Let $\wp \in \spec 
(\pi_{0}A)$ be such that $\chr (k(\wp)) > 0$. Then 
$D_{s}(A|R;k(\wp)) \neq 0$ for infinitely many $s\in {\mathbb N}$ 
provided $A$ fails to be Andr\'e nilpotent at $\wp$.}

\bigskip

\noindent {\bf Remark:} A motivation for the Nilpotence Conjecture 
came from dual topological results centered around
conjectures of Serre and Sullivan as addressed in \cite{Mil,LS,Gro}. See also 
\cite{Tur1} for further speculations on formulating other dual results.

\bigskip

Assuming both of these conjectures, we can now provide a:

\bigskip

\noindent {\it Proof of the Vanishing Conjecture:} First, we have 
$D_{s}(A|R;k(\wp)) = 0$ for $s \gg 0$. Since $\chr (\pi_{0}A) > 0$ 
then $\chr (k(\wp)) > 0$. Thus,
by the Non-Nilpotence Conjecture, $A$ is Andr\'e nilpotent at $\wp$. 
By the Nilpotence Conjecture, $A$ is a homotopy 2-intersection. If, 
additionally, $\fd_{R} (\pi_{*}A) < \infty$ it follows that 
$\pi_{*}(A^{\prime}\otimes_{R^{\prime}}k(\wp))$ is finite and, 
hence, $\Gamma$-nilpotent. It follows from the Nilpotence Conjecture 
that $A$ is a homotopy 1-intersection at $\wp$. \hfill $\Box$

\bigskip 

\section{Proof of the Nilpotence Conjecture at the prime 2}

The goal of this section will be to provide a proof of 
the Nilpotence Conjecture when the base field has characteristic 2. 
This will involve a careful study of a certain map, the character map, 
defined on the homotopy of 
of simplicial supplemented algebras with finite Andr\'e-Quillen 
homology, whose non-triviality implies the Nilpotence Conjecture. In 
fact, in the process of analysing this character map, we will be able 
to establish an upper bound on the top non-trivial degree of the 
Andr\'e-Quillen homology in terms of the non-nilpotence of certain 
operations acting on homotopy. Along the way we will review
some results from \cite{Dwy} and \cite{Goe} and generalize them to 
arbitrary fields of characteristic 2.

\subsection{Connected envelopes and the character map}

We close this section by providing a strategy for proving the 
Nilpotence Conjecture. This will involve first reviewing the 
concept of connected envelopes from \cite{Tur2}. We then 
construct the notion of a character map for connected simplicial 
supplemented algebras with finite Andr\'e-Quillen homology and state 
a conjecture regarding this map whose validity implies the Nilpotence 
Conjecture.

Given $A$ in $s\Ap$, which is connected, we define 
its {\bf connected envelopes} to be a sequence of cofibrations
$$
A = A(1) \stackrel{j_{1}}{\rightarrow} A(2) \stackrel{j_{2}}{\rightarrow} \cdots
\stackrel{j_{n-1}}{\rightarrow} A(n) \stackrel{j_{n}}{\rightarrow} \cdots
$$
with the following properties: 
\begin{itemize}
\item[(1)] For each $n\geq 1$, $A(n)$ is a $(n-1)$-connected. 
\item[(2)] For $s \geq n$,
$$
H^Q_s A(n) \cong H^Q_sA.
$$
\item[(3)] There is a cofibration sequence
$$
S_{\ell}(H^Q_{n} A, n) \stackrel{f_{n}}{\rarrow} A(n) \stackrel{j_{n}}{\rarrow} 
A(n+1).
$$
\end{itemize}
Here we write, for $B \in s\Ap$, $H^{Q}_{s}(A) := D_{s}(A|\ell;\ell)$ 
and, for $V \in \calV_{\ell}$, $S_{\ell}(V,m) := S_{\ell}(K(V,m))$. 
Existence of connected envelopes is proved in \cite[\S 2]{Tur2}.

\bigskip

\noindent {\bf Note:} Paul Goerss has pointed out that connected 
envelopes can also be constructed through a ``reverse'' decomposition 
via collapsing skeleta on the canonical CW approximation.

\bigskip

Now, for $A \in s\Ap$ connected, define the {\it Andr\'e-Quillen 
dimension} of $A$ to be
$$
\aqdim (A) = \mbox{max}\{m\in {\mathbb N}\, | \, H^{Q}_{m}(A)\neq 0 \}.
$$
Assume that $n = \aqdim A < \infty$. Then 
$$
A(n) \simeq S_{\ell}(H^{Q}_{n}(A), n).
$$
Cf. \cite[(2.1.3)]{Tur2}. Summarizing, we have

\begin{proposition}
	For $A \in s\Ap$ connected and $\aqdim A < \infty$ there is a 
	natural map $$\phi_{A} : A \to S_{\ell}(H^{Q}_{n}(A),n),$$ where 
	$n = \aqdim A$, with the property that $H^{Q}_{n}(\phi_{A})$ is an 
	isomorphism. 
\end{proposition}	

Now, assuming $\chr (\ell) > 0$, we noted that $\pi_{*}B$ is naturally a 
divided power algebra. Given a divided power algebra $\Lambda$ in 
characteristic $p$, let $J\subset \Lambda$ to be the divided power 
ideal generated by all decomposables $w_{1}w_{2}\ldots w_{r}$ and 
$\gamma_{p}(z)$ with $w_{1},w_{2},\ldots,w_{r},z\in \Lambda_{\geq 
1}$. Define the {\it $\Gamma$-indecomposables} to be 
$$
Q_{\Gamma}\Lambda = \Lambda/J.
$$
Given $A \in s\Ap$ connected and $n = \aqdim A$ finite, we define 
the {\it character map} of $A$ to be
$$
\Phi_{A} = Q_{\Gamma}\pi_{*}(\phi_{A}) : Q_{\Gamma}\pi_{*}A \to 
Q_{\Gamma}\pi_{*}(S_{\ell}(H^{Q}_{n}(A),n)).
$$

Now, for $B \in s\Ap$, the action of the Andr\'e operation 
$\vartheta$ on $\pi_{*}B$ induces an action on $Q_{\Gamma}\pi_{*}B$ 
by the relation (\ref{anddivrel}) and the fact that $\vartheta$ kills 
decomposables of elements of positive degree. Cf. \cite[(8.9)]{Bou}.

\begin{theorem}\label{charmapthm}
Let $A \in s\Ap$ be connected with $\chr (\ell) = 2$ and $H^{Q}_{*}(A)$ a 
non-trivial finite graded $\ell$-module. Then $\Phi_{A}$ is non-trivial.
\end{theorem}

\bigskip

\noindent {\it Proof of Nilpotence Conjecture at the prime 2:} Let $n = \aqdim B$ 
where $B = A^{\prime\prime}\otimes_{R^{\prime\prime}}\ell$ with $\ell = 
k(\wp)$. By \cor{homsymprop} and \cite[(3.5)]{Goe}, $\vartheta$ acts 
non-nilpotently on every non-trivial element of 
$Q_{\Gamma}\pi_{*}(S_{\ell}(H^{Q}_{n}(B),n))$ if $n \geq 3$. Therefore if 
$\pi_{*}B$ is Andr\'e nilpotent then \thm{charmapthm} implies that $n \leq 2$. 
Thus $B$ is a homotopy 2-intersection by \cite[(2.2)]{Tur2}. 

Since $H^{Q}_{*}(B) \cong H^{Q}_{*}(A(\wp))$, if $A$ is additionally 
$\Gamma$-nilpotent at $\wp$ then $A$ is a homotopy 1-intersection at 
$\wp$, as $\pi_{*}A(\wp)$ is free as a divided power algebra. 
\hfill $\Box$

\bigskip

The goal of the rest of this section will be to provide a proof of 
\thm{charmapthm}.

\bigskip

\noindent {\bf Remark:} If $A$ is a simplicial commutative 
$R$-algebra with finite Noetherian homotopy such that $\chr 
(\pi_{0}A) = 2$, then, for each $\wp \in \spec (\pi_{0}A)$, $\chr 
(k(\wp)) = 2$. Thus \thm{charmapthm} coupled with the Non-Nilpotence 
Conjecture implies (1) of the Vanishing Conjecture when $\chr (\pi_{0}A) = 2$.
Furthermore, as an inspection of the proof of the Vanishing 
Conjecture above shows, (2) of the Vanishing Conjecture follows 
directly from the Nilpotence Conjecture. Thus we have an alternative 
proof of \cite[Theorem B]{Tur3} when $\chr (\pi_{0}A) = 2$.

\subsection{Review of homotopy operations in characteristic 2}

Let $A$ be a simplicial commutative algebra of characteristic 2 (and, 
therefore, a simplicial $\F_{2}$-algebra). Associated to $A$ is
a chain complex, $(C(A), \partial)$, where, for 
each $n \in {\mathbb N}$, we have
$$
C(A)_{n} = A_{n}, \quad \partial = \Sigma_{i=0}^{n}d_{i}: C(A)_{n} 
\to C(A)_{n-1}.
$$
It is standard that we have the identity \cite{Maclane}
$$
\pi_{n}A \cong H_{n}(C(A)).
$$

In \cite{Dwy}, W. Dwyer showed the existence of natural chain maps
$$
\Delta^{k} : (C(V)\otimes C(W))_{i+k} \to C(V\otimes W)_{i} \qquad 
0\leq k \leq i,
$$
where $V$ and $W$ are simplicial $\F_{2}$-modules, having the 
following properties:
\begin{enumerate}
	\item $\Delta^{0}+T\Delta^{0}T = \Delta + \phi_{0}$;
	\item $\Delta^{k}+T\Delta^{k}T = \partial\Delta^{k-1}+\Delta^{k-1}\partial$.
\end{enumerate}
Here $T: C(V)\otimes C(W) \to C(W)\otimes C(V)$ is the twist map, 
$\Delta : C(V)\otimes C(W) \to C(V\otimes W)$ is the shuffle
map \cite[p. 243]{Maclane}, and $\phi_{k} : C(V)\otimes C(W) \to C(V\otimes W)$ is the 
degree ($-k$) map defined by
$$
\phi_{k}(v\otimes w) = 
\begin{cases}
	0 & \mbox{deg}\; v \neq k \; \text{or} \; \mbox{deg}\; w \neq k;\\
	v\otimes w & \text{otherwise}.
\end{cases}
$$

\bigskip

\noindent {\bf Note:} Tensor product of chain complexes is graded 
tensor product and tensor product of simplicial modules is levelwise 
tensor product.

\bigskip

Now, for $x\in C(A)_{n}$ and $1 \leq i \leq n$, define $\Theta_{i}(x) 
\in C(A)_{n+i}$ by $\Theta_{i}(x) = \alpha_{n-i}(x)$ where
$$
\alpha_{t}(x) = \mu\Delta^{t}(x\otimes x) + 
\mu\Delta^{t-1}(x\otimes\partial x),
$$
and
$$
\alpha_{0}(x) = \mu\Delta^{0}(x\otimes x),
$$
where $\mu$ is the map $C(A\otimes A) \to C(A)$ induced by the 
product on $A$. As shown in \cite[\S 3]{Goe}, these natural maps have the 
following properties:
\begin{enumerate}
	\item $\partial \Theta_{i}(x) = \Theta_{i}(\partial x)$ for $2\leq 
	i \leq n$;
	\item $\partial\Theta_{n}(x) = \mu\Delta(x\otimes \partial x)$;
	\item $\partial\Theta_{1}(x) = \Theta_{1}(\partial x) + x^{2}$;
	\item $\Theta_{i}(x+y) = \Theta_{i}(x)+\Theta_{i}(x) +
	      \begin{cases}
			  \partial\mu\Delta^{n-i-1}(x\otimes y) & 2\leq i < n;\\
			  \mu\Delta (x\otimes y) & i = n.
		  \end{cases}$
\end{enumerate} 		  

From these chain properties for the $\Theta_{i}$, there are induced 
homotopy operations
$$
\delta_{i} : \pi_{n}A \to \pi_{n+i}A \qquad 2\leq i \leq n,
$$
or, upon letting $\alpha_{t} = \delta_{n-t}$, we have
$$
\alpha_{t} : \pi_{n}A \to \pi_{2n-t}A \qquad 0\leq t \leq n-2.
$$
Note, in particular, that
\begin{equation}\label{andopid}
	\vartheta = \alpha_{1}.
\end{equation}

The following is proved in \cite{Dwy,GL}:

\begin{theorem}\label{dwythm}
	The homotopy operations $\delta_{i}$ have the following properties:
	\begin{enumerate}
		\item $\delta_{i}$ is a homomorphism for $2\leq i \leq n-1$ and 
		$\delta_{n} = \gamma_{2}$ - the divided square;
		\item $\delta_{i}$ acts on products as follows:
		$$\delta_{i}(xy) = 
		\begin{cases}
			\delta_{i}(x)y^{2} & \mbox{deg}\; y = 0;\\
			x^{2}\delta_{i}(y) & \mbox{deg}\; x = 0;\\
			0 & \text{otherwise};
		\end{cases}$$
		\item if $i < 2j$, then
		$$
		\delta_{i}\delta_{j} = \sum_{\frac{i+1}{2}\leq k\leq 
           \frac{i+j}{3}}\binom{j-i+k-1}{j-s}\delta_{i+j-k}\delta_{k}.
		$$
	\end{enumerate}
\end{theorem}	

\begin{corollary}\label{dwycor}
	The homotopy operations $\alpha_{t}$ have the following properties:
	\begin{enumerate}
		\item $\alpha_{t}$ is a homomorphism for $1\leq i \leq n-2$ and 
		$\alpha_{0} = \gamma_{2}$ - the divided square;
		\item $\alpha_{t}$ acts on products as follows:
		$$\alpha_{t}(xy) = 
		\begin{cases}
			\alpha_{t}(x)y^{2} & \mbox{deg}\; y = 0;\\
			x^{2}\alpha_{t}(y) & \mbox{deg}\; x = 0;\\
			0 & \text{otherwise};
		\end{cases}$$
		\item if $s > t$, then
		$$\alpha_{s}\alpha_{t} = \sum_{\frac{s+2t}{3}\leq q\leq 
             \frac{s+t-1}{2}}\binom{s-q-1}{q-t}\alpha_{s+2t-2q}\alpha_{q}.$$

	\end{enumerate}
\end{corollary}	

\noindent {\it Proof of \cor{dwycor}:} The first two items follow immediately from 
\thm{dwythm} using the identity $\alpha_{t}(x) = \delta_{n-t}(x)$ 
where $\mbox{deg}\; x = n$. The last relation follows from (3) of \thm{dwythm} 
upon letting $j = n-t$, $i = 2n-s-t$, and $k = n-q$. \hfill $\Box$

\bigskip

Our goal at present is to describe homotopy operations for 
simplicial algebras over general fields $\ell$ of characteristic 2. 
Specifically, we will prove:

\begin{theorem}\label{opthm}
	Let $A$ be a simplicial supplemented $\ell$-algebra with $\chr 
	(\ell) = 2$. Then, for $2\leq i \leq n$, the natural operation 
	$\delta_{i} : \pi_{n}A \to \pi_{n+i}A$ satifies properties (1) - (3) 
	of \thm{dwythm}. In particular, for $a,b \in \ell$ and $x,y \in 
	\pi_{n}A$ we have
	$$
	\delta_{i}(ax+by) = a^{2}\delta_{i}(x) + b^{2}\delta_{i}(y) +
	\begin{cases}
		(ab)(xy) & i = n;\\
		0        & otherwise.
	\end{cases}
	$$
	and for $u,v \in \pi_{*}A$
	$$\delta_{i}((au)(bv)) = 
			\begin{cases}
				(ab)^{2}(\delta_{i}(u)v^{2}) & \mbox{deg}\; v = 0;\\
				(ab)^{2}(u^{2}\delta_{i}(v)) & \mbox{deg}\; u = 0;\\
				0 & \text{otherwise}.
			\end{cases}$$
	Furthermore, homotopy operations $\pi_{n}A \to \pi_{n+k}A$, as 
	natural maps of functors of simplicial supplemented $\ell$-algebras, 
	are determined algebraically over $\ell$ by the operations 
	$\delta_{i_{1}}\delta_{i_{2}}\ldots\delta_{i_{r}}$
	with $(i_{1},\ldots,i_{r})$ an admissible sequence of degree k 
	and excess $\leq n$.
\end{theorem}

Recall that the {\it degree} of $I = (i_{1},\ldots,i_{r})$ is 
$i_{1}+\ldots + I_{r}$ and the {\it excess} of $I$ is 
$i_{1}-i_{2}-\ldots -i_{r}$. We will write throughout $\delta_{I} = 
\delta_{i_{1}}\delta_{i_{2}}\ldots\delta_{i_{r}}$. Finally, we call 
$I$ {\it admissible} provided $i_{q-1}\geq 2i_{q}$ for all $2\leq q 
\leq r$.

To prove \thm{opthm}, we need two lemmas. 

\begin{lemma}\label{oplemma1}
	For $n\geq 1$, we have
	$$
	\pi_{*}S_{\F_{2}}(n) \cong \Gamma[\delta_{I}(\iota_{n})|\; 
	\mbox{excess}\, (I) < n].
	$$
	It follows that for any field $\ell$ of characteristic 2
	$$
	\pi_{*}S_{\ell}(n) \cong \Gamma_{\ell}[\delta_{I}(\iota_{n})|\; 
	\mbox{excess}\, (I) < n].
	$$
\end{lemma}

\bigskip

\noindent {\it Proof:} For the first statement, see \cite[\S 7]{Bou} 
or \cite[Remark 2.3]{Dwy}. The second statement follows from the first 
and \prop{symprop}. \hfill $\Box$

\bigskip

For the following, see \cite[12.4.2]{Goe}.

\begin{lemma}\label{oplemma2}
	Let $A$ and $B$ be simplicial commutative $\F_{2}$-algebras. 
	Then the induced action of $\delta_{i}$ on $\pi_{*}(A)\otimes_{\F
	_{2}}\pi_{*}B$ is determined by
	$$\delta_{i}(x\otimes y) = 
	\begin{cases}
		\delta_{i}(x)\otimes y^{2} & \mbox{deg}\; y = 0;\\
		x^{2}\otimes \delta_{i}(y) & \mbox{deg}\; x = 0;\\
		0 & \text{otherwise}.
	\end{cases}$$	
\end{lemma}

\bigskip

\noindent {\it Proof of \thm{opthm}:} Since $\ell$ has characteristic 
2, the operations $\delta_{i}$ are defined on $\pi_{n}A$ and satisfy 
(1) through (3) of \thm{dwythm}. In particular, to compute 
$\delta_{i}(ax+by)$ it is enough, by \cor{homsymprop}, to compute 
$$\delta_{i}(a\iota_{n}\otimes_{\ell} 1 + 1 \otimes_{\ell} b\iota_{n}) \in 
\pi_{*}(S_{\ell}(n))\otimes_{\ell}\pi_{*}(S_{\ell}(n)).$$ Under the 
isomorphism (using \prop{symprop} and Kunneth Theorem) 
$$\pi_{*}(S_{\ell}(n))\otimes_{\ell}\pi_{*}(S_{\ell}(n)) \cong 
(\pi_{*}(S_{\F_{2}}(n)\otimes_{\F_{2}}\pi_{*}(S_{\F_{2}}(n)))\otimes_{\F
_{2}}\ell,$$ $\delta_{i}(a\iota_{n}\otimes_{\ell} 1 + 1 
\otimes_{\ell} b\iota_{n})$ 
corresponds to $\delta_{i}((\iota_{n}\otimes_{\F_{2}} 1)\otimes_{\F_{2}} a + 
(1 \otimes_{\F_{2}} \iota_{n})\otimes_{\F_{2}} b)$. Thus the desired result follows from 
\lma{oplemma2}. Similarly, to compute $\delta_{i}((au)(bv))$ it is 
enough to compute $\delta_{i}((a\iota_{m})\otimes_{\ell}(b\iota_{n})) 
\in \pi_{*}(S_{\ell}(m))\otimes_{\ell}\pi_{*}(S_{\ell}(n))$, or, equivalently, 
$\delta_{i}((\iota_{m}\otimes_{\F2}\iota_{n})\otimes_{\F_{2}} (ab)) 
\in (\pi_{*}(S_{\F_{2}}(m))\otimes_{\F_{2}}\pi_{*}(S_{\F_{2}}(n)))\otimes_{\F_{2}}\ell$.
This again can be computed using \lma{oplemma2}. 

Finally, the last statement follows from \cor{homsymprop} and 
\lma{oplemma1}.
\hfill $\Box$

\bigskip

\noindent {\bf Note:} \thm{opthm} shows that the operations 
$\delta_{i}$ and the relations (1) - (3) of \thm{dwythm} completely 
determine the homotopy operations for simplicial supplemented 
algebras over general fields of characteristic 2. Thus the Galois 
group of $\ell$ over $\F_{2}$ produces no new homotopy operation  
of positive degree nor alters the relations between them. This should 
not be surprising as the same considerations is known to hold 
rationally. See \cite[\S 4]{Qui4}.

\subsection{Quillen's spectral sequence}

We now modify the results of \cite[\S 6]{Goe} to enable to use 
Quillen's fundamental spectral sequence \cite{Qui2, Qui3} over general 
fields of characteristic 2. 

To begin, we need to be more explicit about the functors $S_{\ell}(-)$. 
Let $V$ be an $\ell$-module. For $n\in {\mathbb N}$, define 
$S_{\ell,0}(V) = \ell$ and
$$
S_{\ell,n}(V) = \ell\la v_{1}v_{2}\ldots v_{n}|\; v_{i}\in V\ra. 
$$
Then
$$
S_{\ell}(V) \cong \oplus_{n\in {\mathbb N}} S_{\ell,n}(V).
$$

Next, let $W$ be a non-negatively graded $\ell$-module and define 
\begin{equation}\label{symexp1}
	\calS_{\ell}(W) = \Gamma_{\ell}[\delta_{I}(w)\; |\; w\in W, \;I \, 
	\text{admissible}, \; \mbox{excess}(I) < \mbox{deg}\, w]
\end{equation}
which, by \cor{dwycor}, can be expressed as
\begin{equation}\label{symexp2}
	\calS_{\ell}(W) \cong 	\Gamma_{\ell}[\alpha_{1}^{i_{1}}\alpha_{2}^{i_{2}}\ldots\alpha_{n-2}^{i_{n-2}}(w)
	\; |\; w\in W, \; n = \mbox{deg}\, w,\; i_{1},\ldots, i_{n-2}\in 
	{\mathbb Z}_{+}].
\end{equation}

For $u\in \calS_{\ell}(W)$, we define the {\it weight} of $u$, $\wgt 
(u)$, as follows:
$$
\wgt(u) =
\begin{cases}
	0 & \text{if}\; u\in \ell;\\	
	1 & \text{if}\; u\in W;\\
	\wgt (x) + \wgt (y) & \text{if}\; u = xy;\\
	2\wgt (x) & \text{if}\; u = \delta_{i}(x).
\end{cases}
$$
We then define, for $n \in {\mathbb N}$,
$$
\calS_{\ell, n}(W) = \ell\la u\in \calS_{\ell}(W)| \; \wgt (u) = n\ra.
$$

\begin{proposition}\label{symnprop}
	For a simplicial $\F_{2}$-module $V$ and $n\in {\mathbb N}$ there are a natural 
	isomorphisms $$S_{\ell,n}(V\otimes_{\F_{2}}\ell) \cong 
	S_{\F_{2},n}(V)\otimes_{\F_{2}}\ell$$ and $$\calS_{\ell,n}((\pi_{*}V)\otimes_{\F_{2}}\ell) \cong 
	\calS_{\F_{2},n}(\pi_{*}V)\otimes_{\F_{2}}\ell.$$ As a consequence, if $W$ is a 
	simplicial $\ell$-module then
	$$
	\pi_{*}S_{\ell,n}(W) \cong \calS_{\ell,n}(\pi_{*}W).
	$$
\end{proposition}

\noindent {\it Proof:} The first two statements can be proved just as for 
\prop{symprop}. For the last statement, note that \cite[\S 3]{Goe} 
shows that the isomorphism holds when $\ell = \F_{2}$. Note also that 
a standard argument (e.g. via Postnikov towers) shows that there is a 
simplicial set $X$ and a homotopy equivalence $W \simeq \ell\la X 
\ra$. Thus $$\pi_{*}S_{\ell,n}(W) \cong 
\calS_{\F_{2},n}(\pi_{*}V)\otimes_{\F_{2}}\ell$$ where $V = \F_{2}\la X 
\ra$. Since $\pi_{*}W \cong (\pi_{*}V)\otimes_{\F_{2}}\ell$ it  
follows that
$$
\calS_{\ell,n}(\pi_{*}W) \cong \calS_{\F_{2},n}(\pi_{*}V)\otimes_{\F_{2}}\ell.
$$
\hfill $\Box$ 

\bigskip

We now follow \cite[\S 6]{Goe}.
Let $A$ be a simplicial supplemented $\ell$-algebra and let $IA$ be 
its augmentation ideal. We may assume, using the standard model 
category structure \cite[\S II.3]{Qui1}, that $A$ is almost free, 
i.e. $A_{t} \cong S_{\ell}(V_{t})$ for all $t\geq 1$. Furthermore, 
the composite $V_{t}\subseteq IA_{t} \to QA_{t}$ to the 
indecomposables module is an isomorphism. We now form a decreasing 
filtration of $A$: $$F_{s} = (IA)^{s}.$$ For $A$ almost free, 
$$E^{0}_{s}A = F_{s}/F_{s+1} = (IA)^{s}/(IA)^{s+1} \cong S_{\ell,s}(QA).$$ Applying 
homotopy gives a spectral sequence
\begin{equation}\label{qsseqdef}
	E^{1}_{s,t}A = \pi_{t}E^{0}_{t}A \cong \pi_{t}S_{\ell, s}(QA) 
	\Longrightarrow \pi_{t}A
\end{equation}
with differentials
\begin{equation}\label{qsseqdiff}
	d_{r} : E^{r}_{s,t}A \to E^{r}_{s+r,t-1}.
\end{equation}
This is called {\it Quillen's spectral sequence}.

\begin{theorem}\label{qsseqthm}
	For a simplicial supplemented $\ell$-algebra $A$ there is a spectral 
	sequence of algebras $$E^{1}A_{s,t} = \calS_{\ell,s}(H^{Q}_{*}(A))_{t} 
	\Longrightarrow \pi_{t}A$$ with the following properties:
	\begin{enumerate}
		\item The spectral sequence converges if $\pi_{0}A \cong \ell$. In 
		particular, $E^{r}_{s,t}A = 0$ for $t < s$ for all $r\geq 1$.
		\item For $1\leq r \leq \infty$ there are operations $$\delta_{i} : 
		E^{r}_{s,t}A \to E^{r}_{2s,t+i}A \qquad 2\leq i \leq t$$ of 
		indeterminacy 2r-1 with the following properties:
		\begin{enumerate}
			\item If r = 1, then $\delta_{i}$ coincides with the induced 
			operation $\calS_{\ell,s}(H^{Q}_{*}(A))_{t} \to 
			\calS_{\ell,2s}(H^{Q}_{*}(A))_{t+i}$.
			\item If $x \in E^{r}A$ and $2\leq i < t$ then $\delta_{i}(x)$ 
			survives to $E^{2r}A$ and
			$$
			d_{2r}\delta_{i}(x) = \delta_{i}(d_{r}x) 
			$$
			$$
			d_{r}\delta_{t}(x) = xd_{r}x
			$$
			modulo indeterminacy.
			\item The operations on $E^{r}A$ are induced by the operations on 
			$E^{r-1}A$ and the operations on $E^{\infty}A$ are induced by the 
			operations on $E^{r}A$ for $r< \infty$.
			\item The operations on $E^{\infty}A$ are induced by the operations 
			on $\pi_{*}A$.
			\item Up to indeterminacy, the operations on $E^{r}A$ satisfy the 
			properties of \thm{opthm}.

		\end{enumerate}
		
	\end{enumerate}
	
\end{theorem}

Before we indicate a proof of this omnibus result, a word of 
explanation is needed. First, an element $y\in E^{r}_{s,t}A$ is said 
to be defined up to {\it indeterminacy q} provided $y$ is a coset 
representitive for a particular element of $E^{r}_{s,t}A/B^{q}_{s,t}A$ 
where $$B^{q}_{s,t}A\subseteq E^{r}_{s,t}A \qquad q\geq r$$ is the 
$\ell$-module of elements of $E^{r}_{s,t}A$ which survive to 
$E^{q}_{s,t}A$ but have zero residue class. 

Also, if $A$ is almost free, and hence cofibrant as a simplicial 
supplemented $\ell$-algebra, then
$$
\pi_{*}(QA) \cong H^{Q}_{*}(A).
$$
Cf. \cite[\S 1]{Tur2}. 

\bigskip

\noindent {\it Proof:} First, if $A$ is almost free, we have a
pairing 
$$ 
\pi_{t}(S_{\ell,s}(QA))\otimes \pi_{t^{\prime}}(S_{\ell,s^{\prime}}(QA))
\stackrel{(\mu\Delta)_{*}}{\rarrow} \pi_{t+t^{\prime}}(S_{\ell,s+s^{\prime}}(QA))  
$$ 
which gives a pairing
$$
E^{1}_{s,t}A\otimes 
E^{1}_{s^{\prime},t^{\prime}}A \to E^{1}_{s+s^{\prime},t+t^{\prime}}A
$$ 
and induces an algebra structure on the spectral sequence.

For (1), we simply note that if $A$ is connected then $\calS_{\ell, 
s}(H^{Q}_{*}(A)) = 0$ for $t>s$. Convergence now follows from standard 
convergence theorems. Cf. \cite{Maclane}.

For (2), we have a commutative diagram 
$$
\begin{array}{ccccc}
C((IA)^{s})\otimes_{\F_{2}}C((IA)^{s}) & 
\stackrel{\bar\alpha_{t}}{\longrightarrow} & C((IA)^{s}\otimes_{\F_{2}}(IA)^{s}) 
& \longrightarrow & C((IA)^{s}\otimes_{\ell}(IA)^{s}) \\[1mm]
\hspace*{10pt} \sigma \uparrow  \hspace*{15pt}
&&
&&
\hspace*{0pt}  \downarrow \mu \\[1mm]
C((IA)^{s}) & \stackrel{\alpha_{t}}{\longrightarrow} & C((IA)^{2s}) & = & C((IA)^{2s})
\end{array}	
$$
where $\sigma (u) = u\otimes u$ and $\bar\alpha_{t}(a\otimes b) = 
\Delta^{t}(a\otimes b)+ \Delta^{t-1}(a\otimes\partial b)$. This 
induces a map
$$
\Theta_{i} : (IA)^{s} \to (IA)^{2s},
$$
by again setting $\Theta_{i}(u) = \alpha_{n-i}(u)$ where $n = 
\mbox{deg}\; u$. 

Let $x\in E^{r}_{s,t}A$. Then, modulo $(IA)^{s+1}$,  $x$ is represented by $u \in (IA)^{s}$ 
with the property that $\partial u \in (IA)^{s+r}$. The class of $u$ 
is not unique, but may be altered by adding elements $\partial b \in 
(IA)^{s}$ with $b\in (IA)^{s-r+1}$. 

Define $\delta_{i}(x)\in E^{r}_{2s,t+i}A$ to be the residue class 
of $\Theta_{i}(u)$. Since
$$
\partial\Theta_{i}(u) = \Theta_{i}(\partial u) \in (IA)^{2s+2r} 
\qquad 2\leq i < t,
$$
and
$$
\partial\Theta_{t}(u) = \mu\Delta(u\otimes \partial u) \in (IA)^{2s+r}.
$$
Thus $\delta_{i}(x)$ is defined in $E^{r}_{2s,t+i}A$ and survives to 
$E^{2r}A$ with $d_{2r}\delta_{i}(x) = \delta_{i}(d_{r}x)$ for $2\leq i < t$. 
Also $d_{r}\delta_{t}(x) = xd_{r}x$. This gives us (b).

Now we have a commuting diagram
$$
\begin{array}{ccc}
\pi_{t}((IA)^{s}) & 
\stackrel{(\Theta_{i})_{*}}{\longrightarrow} & \pi_{t+i}((IA)^{2s}) \\[1mm]
\hspace*{10pt} \downarrow  \hspace*{15pt}
&&
\hspace*{0pt}  \downarrow  \\[1mm]
\pi_{t}A & \stackrel{\delta_{i}}{\longrightarrow} & \pi_{t+i}A
\end{array}
$$
and an induced diagram
$$
\begin{array}{ccc}
\pi_{t}((IA)^{s}) & \stackrel{(\Theta_{i})_{*}}{\longrightarrow} & \pi_{t+i}((IA)^{2s}) \\[1mm]
\downarrow \hspace*{10pt}
&&
\hspace*{10pt} \downarrow  \\[1mm]
\pi_{t}((IA)^{s}/(IA)^{s+1}) & \stackrel{}{\longrightarrow} &
\pi_{t+i}((IA)^{2s}/(IA)^{2s+1})\\[1mm]
\hspace*{0pt} \downarrow \cong \hspace*{0pt}
&&
\hspace*{0pt} \cong \downarrow \\[1mm]
\calS_{\ell,s}(H^{Q}_{*}(A))_{t} & \stackrel{\delta_{i}}{\longrightarrow} &
\calS_{\ell,2s}(H^{Q}_{*}(A))_{t+i}
\end{array}
$$
It is now straightforward to check (a), (c), (d), and (e).
\hfill $\Box$

\subsection{Non-triviality of the character map}

We now proceed to prove \thm{charmapthm}. We will in fact prove a more 
general theorem. Specifically:

\begin{theorem}\label{charmapthm2}
	Let $A$ be a simplicial supplemented $\ell$-algebra ($\chr (\ell) = 
	2$) such that $H^{Q}_{*}(A)$ is finite graded as an $\ell$-module. 
	Let $n = \aqdim A$ and assume $n\geq 2$. Then there exists $x\in \pi_{*}A$ and $y\neq 0\in 
	H^{Q}_{n}(A)$ such that under the map $$\pi_{*}\phi_{A} : \pi_{*}A 
	\to \pi_{*}S_{\ell}(H^{Q}_{n}(A),n)$$ we have 
	$$(\pi_{*}\phi_{A})(x) = \alpha_{n-2}^{t}(y)$$ for some $t\geq 1$.
\end{theorem}

\noindent {\it Proof of \thm{charmapthm}:} Assume $n = \aqdim A \geq 
3$. Let $y\in H^{Q}_{n}(A)$, and $x \in \pi_{*}A$ satisfy the properties 
of \thm{charmapthm2}. By \eqn{symexp2}, $\alpha_{n-2}^{t}(y) \neq 0$ 
in $Q_{\Gamma}\pi_{*}(S_{\ell}(H^{Q}_{n}(A),n))$ for all $t\geq 1$. 
We conclude that $\Phi_{A}(x) \neq 0$. 

If $n\leq 2$, then $\Phi_{A}$ is a surjection and, hence, non-trivial.
\hfill $\Box$

\bigskip

Now, in order to prove \thm{charmapthm2} we will need to know 
something about the annihilation properties of homotopy operations. 
Specifically, we will focus on composite operations of the form
$$
\theta(s,t) = \delta_{2^{s}}\delta_{2^{s-1}}\ldots\delta_{2^{t+1}} \quad 
s>t
$$
(where we set $\theta(t+1,t) = \delta_{2^{t+1}}$).

\begin{lemma}\label{annlem}
	Let $i\geq 2$ and $t\geq 1$ be such that $2^{t} < i$. Then 
	$\theta(s,t)\delta_{i} = 0$ for $s\gg t$.
\end{lemma}

\noindent {\it Proof:} 
Write $i = 2^{t-1}+n$ with $n\geq 1$. Note first that an application of 
the relation \thm{dwythm} (3) shows that for any $t\geq 1$, 
$$\delta_{2^{t+1}}\delta_{2^{t}+1}=0=\delta_{2^{t+1}}\delta_{2^{t}+2}.$$ 
We thus assume, by induction, that for any $t$ and $0<j<n$, there exists 
$s\gg t$ such that $$\theta(s,t)\delta_{2^{t}+j}= 0.$$ By another 
application of the relation \thm{dwythm} (3), we have
$$
\delta_{2^{t+1}}\delta_{2^{t}+n} = \sum_{1\leq r\leq 
\frac{n}{3}}\binom{n+r-1}{n-r}\delta_{2^{t+1}+n-r}\delta_{2^{t}+r}.
$$
Notice that, for each such $r$, $2^{t+1}< 2^{t+1}+n-r < 2^{t+1}+n$. Thus, by induction, we 
can find $s\gg t+1$ so that
$$
\theta(s,t+1)\big{(}\sum_{1\leq r\leq 
\frac{n}{3}}\binom{n+r-1}{n-r}\delta_{2^{t+1}+n-r}\delta_{2^{t}+r}\big{)} = 0.
$$
We conclude that $$\theta(s,t)\delta_{2^{t}+n}=\theta(s,t+1)\delta_{2^{t+1}}\delta_{2^{t}+n} = 0.$$ 
\hfill $\Box$

\bigskip

\begin{corollary}\label{anncor}
	Let $I = (i_{1},\ldots, i_{k})$ be an admissible sequence and let 
	$t<k$. Then $\theta(s,t)\delta_{I} = 0$ for $s\gg t$.
\end{corollary}

\noindent {\it Proof:} Since $I$ is admissible, then $$i_{1}\geq 
2i_{2}\geq \ldots\geq 2^{k-1}i_{k}\geq 2^{k} > 2^{t}.$$ Thus, by 
\lma{annlem}, $$\theta(s,t)\delta_{I} = 
(\theta(s,t)\delta_{i_{1}})\delta_{i_{2}}\ldots\delta_{i_{k}} = 0$$ 
for $s\gg t$.
\hfill $\Box$

\bigskip

\begin{proposition}\label{alphacycle}
	Let $A$ be a connected simplicial supplemented $\ell$-algebra, 
	$\chr (\ell) = 2$. Let $y\neq 0$ in $E^{1}_{1,n}A \cong 
	H^{Q}_{n}(A)$, $n\geq 2$. Then there exists $s\geq 1$ such that 
	$\alpha^{s}_{n-2}(y) \in E^{1}_{2^{s},n+2^{s+1}-2}A$ survives to 
	$E^{\infty}A$ (though possibly trivially).
\end{proposition}

\noindent {\it Proof:} Choose $m\geq 1$ and suppose 
$\alpha_{n-2}^{m}(y)$ survives to $E^{r}A$, $r\geq 1$. By 
\thm{qsseqthm} (2) (a), we may assume that $r\geq 2^{m}$. Let 
$w = d_{r}([\alpha_{n-2}^{m}(y)])\in 
E^{r}_{2^{m}+r, n+2^{m+1}-3}A$ by (\ref{qsseqdiff}). By \thm{qsseqthm} (1), 
$w = 0$ 
provided $n+2^{m}-2 \leq r$. Thus if $r\geq n+2^{m}-2$ then the class of
$\alpha^{m}_{n-2}(y)$ survives to $E^{\infty}A$ as all subsequent 
differentials will satisfy the same criterion.

Suppose next that $r < n+2^{m}-2$. Write $n+2^{m}-q = r$ with $n \geq q > 2$. Assume, by 
induction, that if for some $m$ the class of $\alpha^{m}_{n-2}(y)$ survives to 
$E^{n+2^{m}-j}_{2^{m},n+2^{m+1}-2}A$ for $q > j$ then there 
exists $s\gg m$ such that the class of $\alpha^{s}_{n-2}(y)$ survives 
to $E^{\infty}A$. Again, let $w = d_{r}([\alpha_{n-2}^{m}(y)])\in 
E^{r}_{2^{m}+r, n+2^{m+1}-3}A$. Choose $u\in E^{1}_{2^{m}+r, 
n+2^{m+1}-3}A$ to represent the class $w$. By \thm{qsseqthm} and
\prop{symprop}, we have
$$
u =
\begin{cases}
	\sum_{I,l}a_{I,l}\delta_{I}(x_{l}) + \sum_{J}b_{J}z_{J}& r = 
	2^{k}-2^{m},\; k>m;\\
	\sum_{J}b_{J}z_{J} & \text{otherwise}
\end{cases}
$$
where $I = (i_{1},\ldots,i_{k})$ and $J = (j_{1},\ldots,j_{r})$ are 
sequences with $I$ admissible, $a_{I,k},b_{J}\in \ell$, and $z_{J} = 
z_{j_{1}}z_{j_{2}}\ldots z_{j_{r+2^{m}}}$ with 
$x_{k},z_{j_{1}},\ldots,z_{j_{r+2^{m}}} \in H^{Q}_{*}(A)$. 

First assume that $r\neq 2^{k}-2^{m}$. Then 
$d_{r}([\alpha_{n-2}^{m}(y)]) = [u] \in E^{r}_{2^{m}+r, n+2^{m+1}-3}A$
with $u \in E^{1}_{2^{m}+r, n+2^{m+1}-3}A$ decomposable. Note again 
that $\mbox{deg}\; u > 2^{m}$. Thus, by \thm{qsseqthm} (2) (b), 
(c), and (e), $d_{2r}(\delta_{2^{m+1}}[\alpha^{m}_{n-2}(y)]) = 
\delta_{2^{m+1}}d_{r}([\alpha^{m}_{n-2}(y)]) = \delta_{2^{m+1}}[u] = 0$.
Thus $[\alpha_{n-2}^{m+1}y]$ survives to $E^{2r+1}_{2^{m+1}, n+2^{m+2}-2}A$.
Now, let $2r+1 = n+2^{m+1}-j$ and recall that $r =n+2^{m}-q \geq 2^{m}$. Then
\begin{eqnarray*}
	j = (n+2^{m+1})-(2n+2^{m+1}-2q)-1 = 2q-n-1 = q - (n-q) -1 < q.
\end{eqnarray*}
Thus, by induction, there exists $s\gg m$ such that 
$[\alpha_{n-2}^{s}(y)]$ survives to $E^{\infty}$. 

Now assume that $r = 2^{k}-2^{m}$ with $k>m$. 
By definitions of $\alpha_{n-2}$ and $\theta(m,t)$, 
$$\alpha_{n-2}^{m}(y) = \theta(m,0)(y).$$ By \thm{qsseqthm} (2) (b) 
and (c), for $e > m$
$$
d_{2^{e-m}r}([\theta(e,0)y]) = \theta(e,m)d_{r}([\theta(m,0)(y)]) = 
\theta(e,m)w. 
$$
By \thm{qsseqthm} (2) (c) and (e) and \thm{opthm}, $\theta(e,m)w$ is 
represented by
$$
\sum_{I,l}a_{I,l}^{2^{e-m}}\theta(e,m)\delta_{I}(z_{l}) \quad 
\text{modulo indeterminacy}.
$$
Note that $2^{m}< \mbox{deg}\; u$ so we can assume there are no decomposables 
in our choice of representative for $\theta(e,m)w$. As indicated 
above, we have for each 
$I = (i_{1},\ldots,i_{k})$ that $k > m$. Thus, by \cor{anncor}, 
since the sum is finite, there exists $e\gg m$ such that
$$
\theta(e,m)\delta_{I} = 0 \quad \text{for all}\:I.
$$ 
Thus $d_{2^{e-m}r}([\theta(e,0)y]) = 0$ modulo indeterminacy. Therefore
$[\theta(e,0)y]$ survives to \\ $E^{2^{e-m}r+1}_{2^{e}, n+2^{e+1}-2}A$, 
so, by the previous case, there exists $s\gg e$ such that 
$[\alpha_{n-2}^{s}(y)] = [\theta(s,0)y]$ survives to $E^{\infty}$.
\hfill $\Box$

\bigskip

\noindent {\it Proof of \thm{charmapthm2}:}  Choose $y \in 
H^{Q}_{n}(A) \cong E^{1}_{1,n}A$ and choose $s\geq 1$ such that 
$\alpha_{n-2}^{s}(y)\in E^{1}A$ survives to $E^{\infty}A$, which 
exists by \prop{alphacycle}. Under the induced map
$$
E^{r}(\phi_{A}) : E^{r}A \to E^{r}S_{\ell}(H^{Q}_{n}(A),n) 
$$
we have
$$
E^{r}(\phi_{A})([\alpha^{s}_{n-2}(y)]) = [\alpha^{s}_{n-2}(y)] 
$$
for all $\infty \geq r\geq 1$. But, since 
$$
E^{1}S_{\ell}(H^{Q}_{n}(A),n) \cong 
E^{\infty}S_{\ell}(H^{Q}_{n}(A),n),
$$
we can conclude that $E^{\infty}(\phi_{A})([\alpha^{s}_{n-2}(y)]) 
\neq 0$. Thus we can find a nontrivial $x\in \pi_{*}A$ which is represented by 
$\alpha^{s}_{n-2}(y)$ in $E^{\infty}A$ such that 
$(\pi_{*}\phi_{A})(x) = \alpha^{s}_{n-2}(y)$.
\hfill $\Box$

\bigskip

\noindent {\bf Remarks:} 
\begin{enumerate}
	
	\item From the proof of \prop{alphacycle}, an algorithm can be made 
	to determine an $s$ such that $\alpha^{s}_{n-2}(x)$ survives to 
	$E^{\infty}A$ for $x\in H^{Q}_{n}(A)$. Choose $m\geq 1$ such that 
	$\alpha^{m}_{n-2}(x)$ survives to $E^{2^{m}+1}A$, guaranteed by 
	\thm{qsseqthm}.2 (b) and \cor{anncor} (see also \cite[(6.9)]{Goe}). 
	Then, using the procedure in the proof, it can be shown
	that $\alpha^{m+n-2}_{n-2}(x)$ survives to $E^{\infty}A$.
	
	\item Following the philosophy of \cite{Tur1}, the 
     reader can conjecture a dual topological version of \thm{charmapthm2} 
     for nilpotent finite Postnikov towers, using connected covers, which
     would further generalize results of Serre from \cite{Serre} in the 
     spirit of \cite{LS,Gro}.

\end{enumerate}

\section{Proof of the Main Theorem}

We now seek to establish a special case of the Vanishing Conjecture 
as described in the overview. The proof will utilize the validity 
of the Nilpotence Conjecture at the prime 2 while avoiding the need 
to evoke the Non-Nilpotence Conjecture. 

Let $A$ be a simplicial commutative $\F_{2}$-algebra and let 
$(C(A),\partial)$ be the associated chain complex. The following is 
proved in \cite{And2,Dwy}. 

\begin{proposition}\label{dponchains}
	The shuffle map $\Delta : C(A)\otimes_{\F_{2}} C(A) \to 
	C(A\otimes_{\F_{2}} A)$ induces a divided power algebra structure on 
	$C(A)$. Specifically, for each $k\in {\mathbb Z}_{+}$, there is a 
	function $\gamma_{k} : C(A)_{n} \to C(A)_{kn}$ satisfying:
	\begin{enumerate}
		
		\item $\gamma_{0}(x) = 1$ and $\gamma_{1}(x) = x$
		
		\item $\gamma_{h}(x)\gamma_{k}(x) = \binom{h+k}{h}\gamma_{h+k}(x)$
		
		\item $\gamma_{k}(x+y) = \sum_{r+s = k}\gamma_{r}(x)\gamma_{s}(x)$
		
		\item $\gamma_{k}(xy) = 0$ for $k\geq 2$ and $x,y \in C(A)_{\geq 1}$
		
		\item $\gamma_{k}(xy) = x^{k}\gamma_{k}(y)$ for $x\in C(A)_{0}$ 
		and $y\in C(A)_{\geq 2}$
		
		\item $\gamma_{k}(\gamma_{2}(x)) = \gamma_{2k}(x)$
		
		\item $\partial\gamma_{k}(x) = (\partial x)\gamma_{k-1}(x)$
		
		\item $u\in C(A)_{n}$ a cycle then, for $[u] \in \pi_{n}A$, 
		$\delta_{n}([u]) = [\gamma_{2}(u)]$.

	\end{enumerate}	

\end{proposition}

Let $A \to B$ be a map of simplicial commutative $\F_{2}$-algebras and 
$\rho : C(A) \to C(B)$ the induced map of 
chain complexes. Then for $u\in C(A)_{n}$ and all $n > i\geq 0$
$$
\rho (\alpha_{i}(u)) = \alpha_{i}(\rho(u))
$$
where $\alpha_{i} = \Theta_{n-i}$. Recall (\ref{andopid}) that 
$\vartheta = \alpha_{1}$.

\begin{lemma}\label{nilcond}
	Let $A \to B$ be a map of simplicial commutative $\F_{2}$-algebras 
	and suppose $\pi_{s}A = 0$ for $s \gg 0$. Let $u \in C(A)_{n}$, 
	$n \geq 3$, such that $\rho(\partial u) = 0$. Then $\rho (u)$ is a 
	cycle in $C(B)$ and $\vartheta^{r} ([\rho (u)]) = 0$
	in $\pi_{*}(B)$ for $r \gg 0$ provided 
	$\gamma^{r}_{2}(\partial u) = 0$ in $C(A)$ for $r\gg 0$.
\end{lemma}

\bigskip

\noindent {\it Proof:} First, in $C(A)$, we have, by an induction using 
the formulas for $\Theta_{i}$ from \S 2.2, that
$$
\partial \vartheta^{r}(u) = \gamma_{2}^{r}(\partial u).
$$
Since $\gamma_{2}^{r}(\partial u) = 0$ for $r\gg 0$ and $H_{s}(C(A)) = 
0$ for $s\gg 0$, it follows that $\vartheta^{r}(u)$ is a boundary in 
$C(A)$ for $r\gg 0$. We conclude that $\vartheta^{r}([\rho (u)]) = 
[\rho(\vartheta^{r}(u)] = 0$ in $\pi_{*}(B)$.
\hfill $\Box$

\bigskip

\begin{corollary}\label{andnilcond}
	Let $A \to B$ be a level-wise surjection of simplicial commutative 
	$\F_{2}$-algebras such that $\gamma_{2}$ acts locally nilpotently 
	on $(\partial C(A))\cap\ker{\rho}$ and $\pi_{s}A = 0$ for $s \gg 0$. 
	Then $B$ is Andr\'e nilpotent.
\end{corollary}

\noindent {\it Proof:} Given $x \in \pi_{n}B$ with $n \geq 3$, let 
$w \in C(B)_{n}$ be a cycle representitive for $x$ and choose $u\in 
C(A)$ such that $\rho (u) = w$. Then $\rho (\partial u) = 0$ and 
$\gamma_{2}^{r}(\partial u) = 0$ for $r \gg 0$ by assumption. Thus 
$\vartheta^{r} (x) = 0$ for $r \gg 0$ by \lma{nilcond}.
\hfill $\Box$

\bigskip

\noindent {\it Proof of the Main Theorem:} Let $A$ be a simplicial 
commutative $R$-algebra with finite Noetherian homotopy, $R$ a 
Cohen-Macaulay ring of characteristic 2, such that $D_{s}(A|R;-) = 0$ 
for $s \gg 0$, as a functor of $\pi_{0}A$-modules. Note that $A$ has an induced 
simplicial $\F_{2}$-algebra structure. Choose $\wp \in \spec 
(\pi_{0}A)$. Choosing a weak homotopy factorization, 
$(R^{\prime\prime},\fm) \to A^{\prime\prime}$, of $A$ at $\wp$, 
which exists by \prop{final}, then, by \prop{final}, \cite[(3.8.3) 
\& (3.10)]{AFH}, and 
\cite[\S 5]{Mat}, $(R^{\prime\prime},\fm)$ is a Cohen-Macaulay ring 
of depth zero and, hence, locally Artin. Thus, by \prop{final}, we may 
simply assume that $R$ is locally Artinian, that $A$ is a cofibrant 
simplicial commutative $R$ -algebra, and that the unit map $R 
\to \pi_{0}A$ is a surjective local homomorphism. We will now show 
that such $A$ is Andr\'e nilpotent at $\wp$. Note that if $\ell = 
k(\wp)$ then $\chr (\ell) = 2$ since $R$ and $A$ have characteristic 2.

Let $B = A\otimes^{{\bf L}}_{R}\ell = A\otimes_{R}\ell$. Then
\begin{equation}\label{chainmod}
	C(B) \cong C(A)\otimes_{R}\ell \cong C(A)/\fm C(A).
\end{equation}
Thus $\rho : C(A) \to C(B)$ is a surjection and $\ker \rho = \fm C(A)$. 
Since $R$ is locally Artinian, 
\begin{equation}\label{artinnil}
	\fm^{s} = 0 \quad s \gg 0.
\end{equation}
Cf. \cite[2.3]{Mat}. Let $a,b \in \fm$ and let $x,y \in C(A)$ of 
degrees $\geq 2$. By \prop{dponchains} (3) and a straightforward 
induction,
$$
\gamma_{2}^{r}(ax+by) = a^{2^{r}}\gamma_{2}^{r}(x) + 
b^{2^{r}}\gamma_{2}^{r}(y) \quad \text{modulo decomposables}.
$$
Thus, by (\ref{artinnil}) and \prop{dponchains} (4), 
$\gamma_{2}^{r}(ax+by) = 0$ for $r \gg 0$. Hence, by a further 
induction, $\gamma_{2}$ acts locally nilpotently on $\fm C(A)$. Therefore, 
by \cor{andnilcond}, $B$ is Andr\'e nilpotent. We conclude, by the 
validity of the Nilpotence Conjecture at the prime 2, that $A$ is a 
homotopy 2-intersection at $\wp$.
\hfill $\Box$

\end{document}